\documentclass[11pt]{amsart}
\usepackage{fullpage}
\usepackage{amssymb}
\usepackage{epsf}
\usepackage{graphics}
\usepackage{graphicx}
\usepackage{multirow}
\usepackage{amsmath}
\usepackage{enumerate}
\usepackage{epsfig}

\newcommand{\old}[1]{{}}

\newcommand{\id}{\text{Id}}
\newcommand{\su}{\subseteq}
\newcommand{\rest}{\restriction}
\newcommand{\lng}{\langle}
\newcommand{\rng}{\rangle}
\newcommand{\pr}{\operatorname {Pr}}
\newcommand{\Ack}{\operatorname  {Ack}}

\newtheorem{defin}{Definition}[section]
\newenvironment{definition}{\begin{defin} \sl}{\end{defin}}

\newtheorem{theo}[defin]{Theorem}
\newtheorem{fact}[defin]{Fact}

\newtheorem*{abstheorem}{Theorem}

\newenvironment{theorem}{\begin{theo} \sl}{\end{theo}}
\newtheorem{lem}[defin]{Lemma}
\newenvironment{lemma}{\begin{lem} \sl}{\end{lem}}
\newtheorem{coro}[defin]{Corollary}
\newenvironment{corollary}{\begin{coro} \sl}{\end{coro}}
\newtheorem{cl}[defin]{Claim}
\newenvironment{claim}{\begin{cl} \sl}{\end{cl}}
\newtheorem{obs}[defin]{Observation}
\newenvironment{observation}{\begin{obs} \sl}{\end{obs}}

\newtheorem{notat}[defin]{Notation}

\newcommand{\comment}[1]{}

\def\N{{\mathbb N}}

\begin{document}

\title{The threshold for Ackermannian Ramsey numbers}

\author{Menachem Kojman}
\address{Department of Mathematics, Ben Gurion University of the Negev}
\email{kojman@math.bgu.ac.il}
\author{Eran Omri}
\address{Department of Computer Science, Ben Gurion University of the Negev}
\email{omrier@cs.bgu.ac.il}
\date{}
 \begin{abstract} For a function $g:\N\to \N$, the \emph{$g$-regressive Ramsey number} of $k$  is the least $N$ so that 
 \[N\stackrel \min \longrightarrow (k)_g.\]
 
 This symbol means: for every $c:[N]^2\to \N$ that satisfies $c(m,n)\le g(\min\{m,n\})$ there is a \emph{min-homogeneous} $H\su N$ of size $k$, that is,   the color $c(m,n)$ of a pair $\{m,n\}\su H$ depends only on $\min\{m,n\}$.
 
  It is known (\cite{km,ks}) that $\id$-regressive Ramsey numbers   grow  in $k$ as fast as   $\Ack(k)$, Ackermann's function in $k$.  On the other hand, for  constant  $g$,  the $g$-regressive Ramsey numbers   grow  exponentially in $k$, and are therefore primitive recursive in $k$.
 
 We compute below the threshold  in which $g$-regressive Ramsey numbers cease to be primitive recursive and become Ackermannian, by proving:

\begin{abstheorem}
Suppose $g:\N\to \N$ is weakly increasing. Then the $g$-regressive Ramsey numbers are primitive recursive if an only if for every $t>0$ there is some $M_t$ so that for all $n\ge M_t$ it holds that $g(m)<n^{1/t}$ and $M_t$ is bounded by a primitive recursive function in $t$. 
\end{abstheorem}

\end{abstract}

\maketitle

\section{Introduction}
We investigate Ramsey properties of pair-colorings of natural numbers in which the set  of possible colors of a pair depends on the pair. A number $n\in \N$ is identified with the set $\{m\in \N: m<n\}$. 
The set of all two-element subsets of a set $X$ is denoted by $[X]^2$. 

\begin{definition} For a function $f:[\N]^2\to \N$,  an \emph{$f$-coloring of pairs} is a function $c:[\N]^2\to\N$ so that $c(m,n)\le f(m,n)$ for all $\{m,n\}\in [\N]^2$. 
\end{definition}

The standard Ramsey theorems for pairs can be thought of as dealing 
with $f$-colorings for \emph{constant} $f$. When coloring all pairs from $N$ by $C$ colors, there will be a monochromatic subset $B\su N$ of size $k$ if $C^{k\cdot C} \le N$, by the standard proof of Ramsey's theorem.

 On the other hand, if $f$ 
is sufficiently large, that is,  if  $f(m,n)\ge\binom {\max\{m,n\}+1} 2$, then any 
coloring is \emph{equivalent} to an $f$ coloring. Two colorings
 $c_1,c_2$ are equivalent if for all $(m,n), (m',n')$ it holds that
  $c_1(m,n)=c_1(m',n')$ iff $c_2(m,n)=c_2(m',n')$, that is, they induce the same partition of unordered pairs.

The Ramsey behaviour of  colorings of $[\N]^2$ with no limitations at all on the set of colors
 is governed by  the \emph{Canonical Ramsey Theorem} by Erd\H os 
 and Rado, which asserts that for any pair coloring $c:
 [\N]^2 \to \N$ there is an infinite $B\su \N$ so that $c\rest [B]^2$ is 
 \emph{canonical}, that is, is equivalent to one of the following four  colorings: $c_1(m,n)=\min(m,n)$, 
 $c_2(m,n)=\max(m,n)$,   $c_3(m,n)=0$, a   constant coloring or $c_4(m,n)=(m,n)$,  a 1-1 coloring. 
 
 The finite version 
 of the canonical Ramsey theorem asserts that for every $k$ there 
 exists $N$ so that for every $c:[N]^2\to \N$ there is $B\in [N]^k$ so 
 that $c\rest B$ is canonical.    Double exponential upper and lower 
 bounds on $N$ in terms of $k$ are known for the finite canonical Ramsey 
 theorem \cite{lr}. 
 
\bigbreak 

We are interested here in $f$-colorings where $f(m,n)$ depends only 
on $\min\{m,n\}$, that is, when $f(m,n)=g(\min\{m,n\})$ for some function $g:\N\to \N$. When $g=\id$,  such a coloring is called 
\emph{regressive}. In other words, $c$ is regressive if $c(m,n)\le \min\{m,n\}$. More generally, we
 say  that a coloring $c$ is $g$-regressive if $c(m,n)\le g(\min\{m,n\})$.

A set $B\su \N$ is \emph{min-homogeneous} for a coloring $c$ if $c(m,n)$ depends only on $\min\{m,n\}$ for all $(m,n)\in [B]^2$. The important feature of min-homogeneoity is that no matter how large a function $g$ is, a $g$-regressive  
   min-homogeneity Ramsey number exists for every $k$: 
   
\begin{fact} Let $g:\N\to \N$ be arbitrary. Then
\begin{enumerate} 
\item for every $g$-regressive coloring $c:[\N]^2\to \N$  there is an infinite $B\su 
\N$ such that $c\rest [B]^2$ is min-homogeneous.
\item 
 for every  $k$ there is some $N$ so that for every $g$-regressive coloring $c:[N]^2\to \N$ there is a min-homogeneous $B\su N$ of size at least $k$. 
 \end{enumerate}
\end{fact}

\begin{proof} The first item follows from the infinite canonical Ramsey
 theorem, since a regressive coloring cannot  be equivalent neither to 
 $\max\{m,n\}$ nor to a 1-1 coloring on an infinite set.  The second 
 item follows from the first via compactness.   
\end{proof}

Let us introduce the suitable symbolic notation for discussing $g$-Regressive colorings. 

\begin{definition} Let $g:\N\to \N$ be a function. 
Then:
\begin{enumerate}
\item The symbol $N\longrightarrow (k)_g$ means: for every $g$-regressive colorings $c:[N]^2\to \N$ there is a homogeneous $B\su N$ of size $k$. 

\item The symbol $N\overset \min \longrightarrow (k)_g$ means: for every $g$-regressive coloring $c:[N]\to \N$ there is a min-homogeneous $B\su  N$ of size $k$. 
\end{enumerate}
\end{definition}

The $g$-regressive Ramsey theorem (for pairs) is the statement 
\[(\forall k)(\exists N)(N\overset \min \longrightarrow (k)_g )
\]

Recall that the standard proof of Ramsey's theorem gives, for the constant number of colors  $C$, 

\[ C^k\stackrel \min \longrightarrow (k)_C, \qquad C^{k\cdot C} \longrightarrow (k)_C\]

 \bigbreak 
For any function $f:\N\to \N$ the function $f^{(n)}$ is defined by $f^{(0)}(x)=x$ and $f^{(n+1)}(x)=f(f^{(n)}(x))$.  
We recall that Ackermann's function is defined as $Ack(n)=A_n(n)$ where each $A_n$ is the standard \emph{$n$-th approximation} of the Ackermann function, defined by:
 
 \begin{gather}
 A_1(n)=n+1 \\
 A_{i+1}(n)= A_i^{(n)}(n)\notag
 \end{gather}
 
  It is well known (see    e.g. \cite{c}) that each approximation $A_n$
   is \emph{primitive recursive} and that every primitive recursive 
   function is eventually dominated by some $A_n$. Thus Ackermann's 
   function eventually dominates every  primitive recursive function and
    is  *truly* rapidly growing\footnote{See \cite{grs}, p. 60, for an amusing attempt to 
    ``grasp the magnitude" of $A_5(5)$.}. 
   
   Ackermannian lower and upper bounds on $N$ in terms of $k$ are known for the regressive
   Ramsey theorem for  $g=\id$.
 	 This was first proved using methods from mathematical logic and then elementarily \cite{km,ks}. 
 
 \medskip
 
   We are interested here in locating the threshold  for the formidable leap  from a  primitive
   recursive upper bound  to an Ackermannian lower bound in the $g$-regressive Ramsey theorem.
   This threshold obviously lies between the constant functions and $\id$.  
   
   We shall see below that if $g(m)\le m^{\frac{1}{\beta(x)}}$ for some   unbounded and increasing 
     function  $\beta:\N\to \N$ and $\beta^{-1}$ is bounded by a primitive recursive function $f$, then the  $g$-regressive Ramsey numbers  
     are  dominated by  $f$; but if $g(m)=m^{1/\beta(m)}$ where $\beta$ grows to infinity sufficiently slowly,   that is, when $\beta^{-1}$ is Ackermannian, then the $g$-regressive Ramsey number are Ackermannian. 
  
  All functions below are from $\N$ to $\N$ and whenever an expression $x$ may not be integer, it is intended to be replaced by $\lfloor x\rfloor$.
 
\section{The Results}

\subsection{Min-Homogeneity}

For any unbounded function $ \beta: \N \to \N$ define 
 \begin{align} \beta^{-1}(t)= \min(\{n : \beta(n) \ge t\}) \notag \\
\end{align}

%
%
%
%
%
%
%
%
%
%


\begin{notat} For a given $g:\N\to \N$, let $ \nu_g(k)$ denote the
least $N$ which satisfies $ N\stackrel{\min}{\rightarrow}
(k)_{g}$.
\end{notat}

\begin{theo}
 \label{claim:minHom2} 
Suppose $g:\N\to \N$ and $\beta:\N\to \N$ are nonzero,    weakly increasing and
 $g (n)\le   n^{ 1/\beta(n)}$ for all $n$ . Then for every $k\in\N$ it holds 
 that $\nu_{g }(k)\le \beta^{-1}(k).$ 
\end{theo}

\begin{proof}
 Given $1<k \in \N$ let   $N=\beta^{-1}(k)$ and we will show 
 that $N \stackrel \min \longrightarrow (k)_{g}$. Let $C=  g(N) $.
 Since $g$ is increasing, $g (m)\le g (N)$ for all $m < N$. Thus it suffices to show
 that for every coloring $c:[N]\to C$ there exists a min-homogeneous 
 $B\su N$ of size $k$. This holds if $C^k \le N$. Since $C=g(N)\le N^{1/\beta(N)}$ it suffices to show that $(N^{1/\beta(N)})^k\le N$ --- which is obvious, since $\beta(N)\ge k$.   
  \end{proof}

 \begin{corollary}\label{pr}
 Soppose $g$ and $\beta$ are weakly increasing, $g(n)\le n^{1/\beta(n)}$ for all $n$ and $\beta^{-1}$ is bounded by a primitive recursive function. Then $\nu_g$ is bounded by a primitive recursive function. If, furthermore, $g$ is primitive recursive, then $\nu_g$ is primitive recursive.
 \end{corollary}
 
\begin{proof} By the previous theorem $\nu_g$ is bounded by $\beta^{-1}$ and thus is bounded by a primitive recursive function. Since the relation  $ N\stackrel{\min}{\rightarrow}
(k)_{g}$ is primitive recursive when $g$ is, the computation of $\nu_g$ requires a bounded search for a primitive recursive relation and therefore $\nu_g$ is primitive recursive.
\end{proof}

We now begin working towards the proof of the converse of  Corollary \ref{pr}: to show that if $\beta^{-1}$ is Ackermann  and $g(n)=n^{1/\beta(n)}$ then $\nu_g$ is Ackermannian. 
We begin by proving  the special case that   $g(n)=n^{1/\beta(n)}$ and $\beta(n)$ is \emph{bounded}.
%
%
%
%
%



\begin{lemma} 
\label{claim:minHom1} For every $t>0$ let $g_t(n)=n^{1/t}$.  Then the function $\nu_{g_t}$ 
        eventually dominates every primitive recursive function for all $t>0$. 
\end{lemma}

\begin{proof}  The proof is by induction on $t>0$.

The proof involves constructing a ``bad" $f_t$-regressive coloring for $t\ge 1$ by a generalization of the method of construction of a bad $\id$-regressive coloring  in  \cite{ks}. 


\begin{defin}
For a given $t \in \N\setminus \{0\}$, we define a sequence of functions $(f_t)_i : \N \to \N $ as  follows.
 
\begin{align}\label{defEq}
        (f_t)_1(n)&=n+1  \\
        (f_t)_{i+1}(n)&=(f_t)_i^{(\lfloor {n}^{1/t} \rfloor)}(n)
         \end{align}
\end{defin}

\begin{claim} \label{domination}
 For all $0 < t \in \N$ the function $f_t(k)=(f_t)_k(k)$ eventually dominates
 every primitive recursive function.
\end{claim}

\begin {proof}
 By induction on $t$.
 
For $t = 1$ the functions $(f_t)_k=A_k$, the standard   $k$-th
approximations of Ackermann's functions, so every primitive recursive function  is eventually dominated by $f_t(k)$ (see e.g. \cite{c}).

\begin{claim} \label{induction-step}
Let $t > 0$. For all $n > 2^t$, $i > 0$ it holds that 
$(f_{t+1})_{i+2t+2}(n^2) > ((f_t)_i(n))^2 $. 
\end {claim}

We  prove claim \ref{induction-step} by induction on $i$. For $i=1$ we  need  the following:

\begin{observation} \label {obs:prei1}
For every $t,k,n > 0$ it holds that $(f_{t})_{k}(n)  \ge  n + (\lfloor {n}^{1/t} \rfloor)^{k-1}$
\end{observation}

\begin{proof}
 We show observation \ref{obs:prei1} by induction on $k$. If $k = 1$, it follows by definition
  that $(f_{t})_{k}(n)  = n+1 =  n + (\lfloor {n}^{1/t} \rfloor)^{k-1}$. Let $k>1$.
 By definition $(f_{t})_{k+1}(n)  = (f_t)_k^{(\lfloor {n}^{1/t} \rfloor)}(n)$ 
 and by applying the induction hypothesis $\lfloor {n}^{1/t} \rfloor$ times
 we get that the right hand side of the equation is larger than 
 $n+ ((\lfloor {n}^{1/t} \rfloor) (\lfloor {n}^{1/t} \rfloor)^{k-1})$ which is
 $n+ (\lfloor {n}^{1/t} \rfloor)^{k}$.
\end{proof}

\begin{observation} \label {obs:i1}
$(f_{t+1})_{2t+3}(n^2)  >  n^2 + 2n + 1$
\end{observation}

\begin{proof}
 By observation \ref{obs:prei1} we have that $(f_{t+1})_{2t+3}(n^2)  \ge n^2 + (\lfloor {n}^\frac{2}{t+1} \rfloor)^{2t+2}$. 
 
 Now
\begin{multline} \notag n^2 + (\lfloor {n}^\frac{2}{t+1} \rfloor)^{2t+2} \ge
n^2 + ({n}^\frac{2}{t+1} - 1)^{2(t+1)} \cr \ge n^2 +  (n^\frac {4}{t+1} - 2n^\frac{2}{t+1} + 1)^{t+1} \cr > n^2 +  (n^\frac {2}{t+1}( n^\frac{2}{t+1} - 2))^{t+1} > 2n^2 > n^2 + 2n + 1
\end{multline}
\end{proof}

When $i = 1$, by observation \ref{obs:i1},  $(f_{t+1})_{i+2t+2}(n^2) = (f_{t+1})_{2t+3}(n^2) > n^2 + 2n + 1 = {((f_{t})_1(n))}^2 = {((f_{t})_i(n))}^2$

We now assume that claim \ref{induction-step} is true for $i$ and prove it for $i+1$.
\begin {claim} \label{monotonicity}
$ \forall j \in \N$  $(f_{t+1})_{i+2t+2}^{(j)}(n^2) > ((f_t)_i^{(j)}(n))^2 $.
\end{claim}

\begin{proof}
We show claim \ref{monotonicity} by induction on $j$. For $j = 1$ the claim is induced by 
the induction hypothesis for $i$. For $j > 1$ we have
$(f_{t+1})_{i+2t+2}^{(j+1)}(n^2)$ $=$ $(f_{t+1})_{i+2t+2}((f_{t+1})_{i+2t+2}^{(j)}(n^2))$. 
The latter term is larger than $(f_{t+1})_{i+2t+2} (((f_t)_i^{(j)}(n))^2)$ by monotonicity and the 
induction hypothesis for $j$. Now, if we denote $n' = (f_t)_i^{(j)}(n)$, we easily see, by the induction hypothesis for j or for i, that  
$(f_{t+1})_{i+2t+2} (((f_t)_i^{(j)}(n))^2)$ $>$ $((f_t)_i ((f_t)_i^{(j)}(n)))^2$ which is, 
in fact, $((f_t)_i^{(j+1)}(n))^2$
\end{proof}

We still need to show the induction step for claim \ref{induction-step}
$(f_{t+1})_{i+1+2t+2}(n^2)$ $=$ $(f_{t+1})_{i+2t+2}^{(\lfloor {n}^\frac{2}{t+1} \rfloor)}(n^2)$ 
 $\ge$ $(f_{t+1})_{i+2t+2}^{(\lfloor {n}^{1/t} \rfloor)}(n^2)$ and by claim \ref{monotonicity},
  the latter term is larger than $((f_t)_{i}^{(\lfloor {n}^{1/t} \rfloor)}(n))^2 $ $=$ $((f_t)_{i+1}(n))^2 $

 That concludes the proof of claim  \ref{induction-step} and therefore also of claim \ref{domination}.
 \end{proof}

We turn now to the construction of bad $g_t$-regressive colorings. 

 For a given natural number $k > 2$ and a given 
 $g:\N \to \N$ that is monotonically increasing such that for some $t \in \N$ it 
 holds that $k \le \lfloor  \frac {\sqrt {g(t)}} {2} \rfloor$, we define a 
 sequence of functions  $(f_g)_i : \N \to \N$ as follows.

 \begin{defin} \label{defCol}
 Let $\mu = \mu_g(k) = \min(\{t \in \N : k \le \lfloor  \frac {\sqrt {g(t)}} {2} \rfloor \})$

 and let
 \begin{align}       
 (f_g)_1(n) = n+1 \\
 (f_g)_{i+1}(n) = (f_g)_i^{(\lfloor \frac {\sqrt {g(n)}} {2} \rfloor)}(n)
 \end{align}
 \end {defin}

Define a sequence of semi-metrics $\lng {(d_g)}_i : i\in \N\rng$ on 
$\{n : n \ge \mu \}$ by setting, for $m,n \ge \mu$, 

\begin{gather}
{(d_g)}_i(m,n) = |\{l \in \N : m < (f_g)_i^{(l)} (\mu) \le n\}|
\end{gather}

 For $n > m \ge \mu$ let $I_g(m,n)$ be the greatest $i$ for which ${(d_g)}_i(m,n)$ 
is positive, and  $D_g(m,n)={(d_g)}_{I(m,n)}(m,n)$.

Let us fix the following (standard) pairing function $\pr$ on 
$\N^2$ $$\Pr(m,n)=\displaystyle \binom{m+n+1}{2}+n$$

$\pr$ is a bijection between $[\N]^2$ and $\N$ and is monotone in each 
variable.  Observe that if $m,n\le l$ then $ \Pr(m,n)<4l^2$ for all 
$l>2$.

Define a pair coloring $c_g$ on $\{n:n\ge \mu \}$ as follows: 

\begin{gather}
\label{coloring}
c_g(\{m,n\}) = \Pr(I_g(m,n),D_g(m,n)) 
\end{gather}

\begin{claim} \label{smallDg}  
For all $n > m \ge \mu$, $D_g(m,n) \le \frac {\sqrt {g(m)}} {2} $.
\end{claim}

\begin{proof}
  Let $i=I_g(m,n)$. Since ${(d_g)}_{i+1}(m,n)=0$, there exist $t$ 
 and $l$ such that $t = (f_g)_{i+1}^{(l)}(\mu) \le m < 
 n < (f_g)_{i+1}^{(l+1)}(\mu) = (f_g)_{i+1}(t)$.  
 But $(f_g)_{i+1}(t) = (f_g)_{i}^{(\lfloor \frac {\sqrt {g(t)}} {2} \rfloor)}(t)$ 
 and therefore $\frac {\sqrt {g(t)}} {2}  \ge {(d_g)}_i(t,(f_g)_{i+1}(t))\ge D_g(m,n)$.
\end{proof}

\begin{claim} \label{g_regressive} 
$c_g$  is $g$-regressive on the interval $[\mu ,(f_g)_k(\mu))$.
\end{claim}

\begin{proof} Clearly, ${(d_g)}_k(m,n)=0$ for $\mu \le 
m < n < (f_g)_k(\mu)$ and therefore $I_g(m,n) < k \le \frac {\sqrt {g(m)}} {2} $.  From claim \ref{smallDg} we know that $D_g(m,n)\le \frac {\sqrt {g(m)}} {2}$.  Thus, 
$c_g(\{m,n\}) \le \pr(\lfloor \frac {\sqrt {g(m)}} {2} \rfloor,\lfloor\frac {\sqrt {g(m)}} {2}\rfloor)$, which is 
$< g(m)$, since $\frac {\sqrt {g(m)}}{2} > 2$.
\end{proof}

\begin{claim} \label{noMinHom} 
For every $i \in\ N$, every sequence $x_0 < x_1< \dots <x_i$ that satisfies 
${(d_g)}_i(x_0,x_i)=0$ is not min-homogeneous for $c_g$.
\end{claim}

\begin{proof}
The claim is proved by induction on $i$.  If $i=1$ then there are no 
$x_0 < x_1$ with ${(d_g)}_1(x_0,x_1) = 0$ at all. Let $i > 1$ and suppose 
to the contrary that $x_0 <x_1 < \dots < x_i$ form a 
min-homogeneous sequence with respect to $c_g$ and that ${(d_g)}_i(x_0,x_i) = 0$.  
Necessarily, $I_g(x_0,x_i) = j < i$. By min-homogeneity, 
$I(x_0,x_1)=j$ as well, and ${(d_g)}_j(x_0,x_i)={(d_g)}_j(x_0,x_1)$.  Hence, 
$\{x_1,x_2,\dots x_i\}$ is min-homogeneous with ${(d_g)}_j(x_1,x_i)=0$ --- 
contrary to the induction hypothesis.
\end{proof}

\begin{coro}
 There exists no $H \su [\mu ,(f_g)_k(\mu))$ of size $k+1$ that is min-homogeneous for $c_g$.
\end{coro}

\begin{coro} \label{lowerBoundDomination}
 If the function ${(f_g)}_k(k)$ dominates every primitive recursive function 
 (Ackermannian in terms of $k$) and $\mu_g(k)$ is bound by some primitive 
 recursive function, then the lower bound for min-homogeneity for 
 $g$-regressive colorings also dominates every primitive recursive function.
\end{coro}
 
 \begin{proof}
 The collection of primitive recursive functions is closed under composition. Thus,
 the function ${(f_g)}_k(\mu_g(k)) - \mu_g(k)$ is Ackermannian in terms of $k$.
 Moreover, it is Ackermannian in terms of $\mu_g(k) + k + 1$. Therefore, we may
 allow ourselves to set the color of every pair $(m,n)$ such that $m < \mu_g(k)$ to be 
 $0$ and by that present a $g$-regressive coloring of $[{(f_g)}_k(\mu_g(k))]^2$ 
 that yields no min-homogeneous $H \su [{(f_g)}_k(\mu_g(k))]$ of size 
 ${\mu_g(k) + k + 1}$.  
 \end{proof} 
 
 Now, to conclude the proof of theorem \ref{claim:minHom1} we need only observe
 that for a given $j \in \N$ the function $\frac {k^{\frac {1}{2j}}}{2}$ grows 
 asymptotically faster than $k^{\frac {1}{4j}}$ and therefore, 
 by claim \ref{domination}, for any $j \in \N$ ${(f_g)}_k(k)$  for $g(m) = m^j$
 dominates every primitive recursive function. On the other hand, for such $g$,
 $\mu_g(k) \le 4^j k^{2j}$. Hence, by corollary \ref{lowerBoundDomination} we
 establish that the lower bound for min-homogeneity for 
 $g$-regressive colorings for $g(m) = m^j$ dominates every primitive recursive 
 function.

\end{proof}

%

\begin{theo}
 Let  $\beta^{-1}(n):=\Ack(n+3)$ (so   $\beta$ is basically $\Ack^{-1}$) and let $g(n)=n^{1/\beta(n)}$. There exists a $g$-regressive coloring $c: [\N]^2 \to \N$ 
 such that for every primitive recursive function $f: \N \to \N$ there exists $N_f \in \N$ 
 such that for all $m > N_f$ and  $H \su  m $    which is min-homogeneous 
 for $c$ it holds  that $f(|H|)< m$. 
\end{theo}

\begin {proof}
	We define two increasing sequences $\{k_t\}$ and $\{\mu_t\}$  and then let $\beta(n)=t+1$ if $\mu_t\le n<\mu_{t+1}$. Using the definition of $(f_g)_i$
	given in \ref{defCol}, we define a $g$-regressive coloring $c$, where $g(n)=n^{1/\beta(n)}$, so that in the interval $[\mu_t,\mu_{t+1})$ there is no min-homogeneous set of size $k_t$. 
	 
	We  denote   $g_t(n): = n^{ 1/t} $. Let:
  \begin{align}
  			\mu_0 &= 0  \notag \\
  			\mu_1 &= 10^4 \notag \\
  			k_1 &= 18 \notag
\end{align}

And for all $t>1$,
\begin{align}
  			k_{t} &=   \frac{\sqrt{g_{t-1}(\mu_{t-1})}}{2}   = \frac { \mu_{t-1} ^{1/2(t-1)}} 2\notag\\
         \mu_{t}  &=\Ack(t+3)  
         \notag
  \end{align}
%
%
On $[0,\mu_1)$ we define $c(m,n)$ as follows: color all $\{m,n\}$ from $[0,43)$ regressively by the colors $\{0,1\}$ with no min-homgeneous set of size $12$. This is possible, since the (usual) Ramsey number of $5$ is $\ge 43$, so there is  a 2-coloring of $[1,43)$ with no homogeneous set of size $5$, hence with no min-homogeneous set of size $11$. For $m,n\ge 43$ color as follows: write out $|n-m|$ in base $10$ and let $c(m,n)=\Pr(d_1,d_2+1)$ where $d_1\in \{0,1,2,3\}$ is the maximal power of $10$ smaller than $|n-m|$ and $d_2$ is the $f_1$-th decimal digit. This coloring allows no min-homogeneous sets of size $6$ in $[42,10^4)$. So letting $c(m,n)=0$ for $m< 42$ and $42\le m\le 10^4$, we get that below $\mu_1$ there are no min-homogeneous sets of size $k_1=18$.

Now we need to define $c$ on $[\mu_{t-1},\mu_{t})$ for all $t>1$.
Let $k_{t}=\frac { \mu_{t-1} ^{1/2(t-1)}} 2$.
 Observe that  we may color pairs over 
  the interval $[\mu_{t-1}, (f_{g_{t-1}})_{k_t}(\mu_{t-1}))$ if $\mu \ge 4(g_t(k))^2$ using $c_{g_t}$ 
  with $\mu_t$ instead of $\mu_{g_t}(k)$. This coloring is   $g_t$-regressive with
  no min-homogeneous $H \su [\mu, ((f_{g_t})_k(\mu))$ of size $k+1$. This is true since the proofs of
  claims \ref{g_regressive} and \ref{noMinHom} made no use of the minimality of $\mu_{g_t}(k)$. 
  
  To define  $c$ on $[\mu_{t-1},\mu_{t})$ it suffices, then, to prove:
  
  \begin{claim} \label{intervalGrowth:claim}
  $\Ack(t+3)< (f_{g_{t-1}})_{k_{t}}(\mu_{t-1}))$ for all $t>1$.
  \end{claim}
  
  \begin {proof}
  We first prove claim \ref{intervalGrowth:claim} 
  for $t = 2$. We have, by claim \ref{muGrowth:claim}, that 
  $(f_{g_1})_{k_2}(\mu_1)) > 
	  (A_{k_{2}-18}(\left\lfloor{\mu_{1}}^{\frac{1}{8}}\right\rfloor))^{8}$.
   Since $k_2 = 50$, the latter term is
   $(A_{32}(\left\lfloor{10^4}^{\frac{1}{8}}\right\rfloor))^{8}$ and thus
   larger than $A_{32}(3) > A_{5}(5)$.
     
   Let $t >2$. We know that $\mu_{t-1} = \Ack(t+2)$ 
   and hence it clearly holds that $k_t - 16t^2 + 28t -10 =  
   \frac { \mu_{t-1} ^{1/2(t-1)}} 2 - 16t^2 + 28t -10 = 
   \frac { (\Ack (t+2)) ^{1/2(t-1)}} 2 - 16t^2 + 28t -10 >
   t+3$ and it also clearly holds that 
   $\mu_{t-1} ^ {\frac{1}{2^{4t-5}}} > t+3$. Thus,  by claim 
   \ref{muGrowth:claim}, we have that
   $(f_{g_{t-1}})_{k_{t}}(\mu_{t-1})) \ge 
			(A_{k_{t}-16t^2 + 28t -10}(\left\lfloor{\mu_{t-1}}^{\frac{1}{2^{4t-5}}}\right\rfloor))^{2^{4t-5}} > A_{t+3}(t+3)$.

  \end{proof}

  \begin {claim} \label{muGrowth:claim}
  		For all $t > 1$ it holds that 
  			$(f_{g_{t-1}})_{k_{t}}(\mu_{t-1})) > 
			(A_{k_{t}-16t^2+28t-10}(\left\lfloor{\mu_{t-1}}^{\frac{1}{2^{4t-5}}}\right\rfloor))^{2^{4t-5}}$	
  \end{claim}
  
  \begin {proof}
	Observe that $(A_{k_{t}-16t^2+28t-10}
	(\left\lfloor{\mu_{t-1}}^{\frac{1}{2^{4t-5}}}\right\rfloor))^{2^{4t-5}}$
	is actually $((f_1)_{k_{t}-16t^2+28t-10}
	(\left\lfloor{\mu_{t-1}}^{\frac{1}{2^{4t-5}}}\right\rfloor))^{2^{4t-5}}$
	 Now, by applying claim \ref{induction-step} to the latter term,
	 we get $((f_1)_{k_{t}-16t^2+28t-10}
	(\left\lfloor{\mu_{t-1}}^{\frac{1}{2^{4t-5}}}\right\rfloor))^{2^{4t-5}}
	 < ((f_2)_{k_{t}-16t^2+28t-10+2+2}
	(\left\lfloor{\mu_{t-1}}^{\frac{1}{2^{4t-6}}}\right\rfloor))^{2^{4t-6}}$,
	since the parameter $t$ of claim \ref{induction-step} is $1$
	here. If we apply it now to right hand side term, the parameter
	$t$ of the claim would be $2$ and we would find that the
	latter term is smaller than  
	$((f_3)_{k_{t}-16t^2+28t-10+2+2+4+2}
	(\left\lfloor{\mu_{t-1}}^{\frac{1}{2^{4t-7}}}\right\rfloor))^{2^{4t-7}}$.
	Generally, if we apply the claim $j$ times we get that
	$((f_1)_{k_{t}-16t^2+28t-10}
	(\left\lfloor{\mu_{t-1}}^{\frac{1}{2^{4t-5}}}\right\rfloor))^{2^{4t-5}}
	 < ((f_{j+1})_{k_{t}-16t^2+28t-10+j^2+3j}
	(\left\lfloor{\mu_{t-1}}^{\frac{1}{2^{4t-5-j}}}\right\rfloor))^{2^{4t-5-j}}$
	since we may replace $\sum_{l=1}^{j}2j$ with $j^2+j$.
	Now, if we let $j=4t-5$, we get
	$((f_1)_{k_{t}-16t^2+28t-10}
	(\left\lfloor{\mu_{t-1}}^{\frac{1}{2^{4t-5}}}\right\rfloor))^{2^{4t-5}}
	 < (f_{4(t-1)})_{k_{t}}({\mu_{t-1}})$. Note that we are allowed to 
	 apply  claim \ref{induction-step} $4t-5$ times, only if, for all
	 $1 \le j \le 4t-5$ it holds that 
	 $\mu_{t-1}^{\frac{1}{2^{4t-4-j}}}>2^{j}$,
	 or that $\mu_{t-1} > 2^{j 2^{4t-4-j}}$ and that is true 
	 for all $t>2$ since $\mu_{t-1}$ is clearly larger than
	 $2^{(4{(t-1))2^{4{(t-1)}}}}$. For $t=2$ it is also true and may be
	 easily verified by hand.
	

  	On the other hand, it holds that $(f_{g_{t-1}})_{k_{t}}(\mu_{t-1}) >
	(f_{4(t-1)})_{k_{t}}(\mu_{t-1})$
  	since $\mu_{t-1}$ is larger than $2^{4t}$ for all $t>1$ and
	therefore $\frac{n^{\frac{1}{2(t-1)}}}{2} \ge n^{\frac{1}{4(t-1)}}$.

  \end {proof}

	\begin{observation}
	  The coloring $c$ is $g$-regressive.
	\end{observation}   
	\begin {proof}
		For any $m,n$ such that $\beta(m) = \beta(n) = t$ we know that
		$c(m,n) \le \left\lfloor m^{1/t}\right\rfloor$ since $c(m,n) = c_{g_t}(m,n)$ and $c_{g_t}$ 
		is $g_{t}$-regressive on the interval. Otherwise, $c(m,n) = 0$ which is
		always smaller than $m^{1/\beta(m)}$.
	\end {proof}
	
	\begin{observation}
	  For any given $N \in \N$ with $\beta(N) = j$, there is no min-homogeneous $H \su [N]$ of size 
	  $(k_{j-1}+1)^{2}+18$.
	\end{observation}   
	
	\begin {proof}
	  From claim \ref{muGrowth:claim} it is clear that for all $t > 1$ it holds that $k_t < k_{t+1}$
	  and that $k_t > t$. Thus, since at each interval $[\mu_t,\mu_{t+1})$ for ant $t < j$ there exist 
	  no min-homogeneous subset of size $k_{t}+1$ and hence, no min-homogeneous subset of size $k_{j-1}+1$. 
	  Therefore, in the union of all those intervals there is no  min-homogeneous subset of size 
	  $(k_{t}+1)t < (k_{j-1}+1)^{2}$. Now,   in the first interval  there 
	  can be no no min-homogeneous of size $18$, there is no min-homogeneous $H \su [N]$ of size $(k_{j-1}+1)^{2}+18$ in the union of the first
	  $j$ intervals of which $[N]$ is a subset. 
	\end {proof}
  
  To conclude the proof we only need to observe that given a primitive recursive function $f$, there
  exists a $k_f \in \N$ such that for every $n > 4$ it holds that 
  $A_{k_f}(\left\lfloor n^{1/\lg \lg{n}}\right\rfloor) > f({(n+1)^2+18})$. Now, because
  $k_t$ grows extremely faster than $t$, we can find a $t$ such that $k_{t}-32(t-1)^2-4(t-1)+3 > k_f$ and
  $\lg \lg \mu_t > t$.
  Set $N_f$ to be $\mu_{t+1}$. Given $n > N_f$ with $\beta(n) = j$. We have that $j > t+1$. Assume to the 
  contrary that there exists a min-homogeneous $H \su [n]$ of size $f^{-1}(n)$ then 
  $f^{-1}(n) < (k_{j-1}+1)^{2}+18$. Thus, $n' = \sqrt{f^{-1}(n)-18}-1 < k_{j-1} < \mu_{j-2}$.
  Now, 
  $n = f({(n'+1)^2+ 18}) $ 
  $<  A_{k_f}(\left\lfloor {n'}^{\frac{1}{\lg \lg{n'}}}\right\rfloor)<$
  $ A_{k_{j-1}-32(j-2)^2-4(j-2)+3}(\left\lfloor \mu_{j-2}^\frac{1}{\lg \lg{\mu_{j-2}}}\right\rfloor) $
  $< \mu_{j-1}$. Contrary to $\beta(n) = j$.

\end {proof}

We can now prove the main theorem of the paper:

\begin{theorem}\label{ackcol}
Suppose $g:\N\to \N$ is eventually smaller than $n^{1/t}$ for every constant $t>1$. Then $\nu_g$ is bouded by a primitive recursive function if an only if the least number $M_t$ which satisfies $g(n)<n^{1/t}$ for all $n\ge M_t$  is bounded by a primitive recursive function in $t$. 
\end{theorem}

\begin{proof}Suppose first that $M_t$ is bounded by some primitive recursive function in $t$. Replacing $g(n)$ by $\max\{g(n):m\le n\}$ we may assume that $g$ is weakly increasing and $M_t $ would still be bounded by a primitive recursive function. Now apply Corollary \ref{pr}. This takes care of the ``if" part. 

 The ``only if" part follows directly from Theorem \ref{ackcol} above.
\end{proof}
\subsection{Homogeneity} \label{proof:regHom}

We look now at the threshold $g$ at which one can guarantee the usual Ramsey theorem for $g$-regressive colorings, that is, have homogeneous rather than just min-homogeneous sets.

\begin{theo}  
  \label{claim:regHom1} Suppose  $f:\N\to \N^+$ satisfies 
 ${\displaystyle \lim_{n\to\infty} (f(n)) = \infty}$ 
 and let $g(x)=\frac{\lg x }{f(x) \lg  \lg x}$ for $x\ge 4$ 
 and $g(x)=0$ for $x<4$. Then  $\forall k $ $\exists N$ $ N\to (k)_{g}^2$.
\end{theo}

\begin{proof}
Let $f:\N \to \N^+$  be a function such that 
${\displaystyle \lim_{n\to\infty} (f(n)) = \infty}$, 
and $g(x)=\frac{\lg x }{f(x) \lg  \lg x}$ for $x\ge 4$ and $g(x)=0$ for $x<4$.
Given $k \in \N$, find $N\ge k$ so that $f(N)>k$ and
 $f(N)>f(m)$ for all $m<N$. Such $N$ exists, since ${\displaystyle
 \lim_{m\to\infty} (f(m)) = \infty}$. Since $f(N)>f(m)$ for
  all $m<N$, it follows that $g(m)\le g(N)$ for all $m<N$ as well. So given a
  $g$-regressive coloring
$c : [N]^2\to \N$ we have that $c(m,n)\le g(N)$ for all $(m,n)\in [N]^2$. Put
$C= \lfloor g(N)+1 \rfloor$. If $C=1$ then $N$ itself is homogeneous of size $\ge k$, so assume
that $C\ge 2$. The standard proof of Ramsey's theorem with $C$ colors gives a
homogeneous
 $B\su N$ of size $k$ in case $N> C^{k\cdot C}$, which holds here, since $\lg
 g(N)< \lg \lg N$ and therefore
 
\[ C^{k\cdot C}=2^{\lg C \cdot k\cdot C}=2^{\lg g(N) \cdot k \cdot g(N)} <
2^{\lg \lg N \cdot k \cdot g(N)}\le 2^{\lg N}=N
\]

It should be noted that this is of interest when $f$ grows slowly (e.g.
$f(m)=\log^* m$).
\end{proof}

\begin{theo}
\label{claim:regHom2} For every $s \in \N$ and for $g(i)={\frac{\lg i}{s}}$ it
holds that $\exists k$ $\forall N$ $N \nrightarrow {(k)_{g}^2}$
\end{theo}

\begin{proof}

Let $s \in \N$ and $g(i)={\frac{\lg(i)}{s}}$. We set $k = 2s + 1$
and we show a $g$-regressive
 coloring $C: {\N}^2 \rightarrow \N$ where there exists no $S \subseteq \N$ of
 size $\geq k$
  that is homogeneous for $C$. For any $n \in \N$, let $r_s(n)$ be the
  representation of $n$ in $s$
   basis. For any $m,n \in \N$ such that $m < n$ and
   $\left\lfloor{\log}_s(m)\right\rfloor = \left\lfloor {\log}_s(n)
   \right\rfloor$, let $f(m,n)$ be the smallest index $i$ such that $r_s(m)[i] \neq
   r_s(n)[i]$. We define $c$ as

\[ c(m,n) = \left\{ \begin{array}{ll}
\left\lfloor{\log}_s(m)\right\rfloor & \mbox{if
  $\left\lfloor{\log}_s(m)\right\rfloor  \neq \left\lfloor {\log}_s(n)
  \right\rfloor$};\\
f(m,n) & \mbox{if $\left\lfloor{\log}_s(m)\right\rfloor = \left\lfloor
    {\log}_s(n) \right\rfloor$}.\end{array} \right. \]

\begin{observation} \label{observ:svalues}
Let $Y = \{y_1,y_2,...,y_{s+1}\}$ where $ y_1 < y_2 <...< y_{s+1}$, be a homogeneous
set for $C$. Then
$\left\lfloor{\log}_s(y_1)\right\rfloor  < \left\lfloor {\log}_s(y_{s+1})
\right\rfloor$.
\end{observation}

 To show Observation \ref{observ:svalues}, let $Y$ be a homogeneous
 set for $C$ and suppose to the contrary that 
 $\left\lfloor{\log}_s(y_1)\right\rfloor = \left\lfloor {\log}_s(y_{s+1})
 \right\rfloor$,
  from the definition of $c$ we get that $f$ is constant on $Y$. Thus elements
  of $Y$,
   pairwise differ in the $i$'th value in their $s$ basis representation for
   some index $i$,
    which is impossible since there are only $s$ possible values
    for any index. Contradiction.

Now, Let $X = \{x_1,x_2,...,x_{2s+1}\}$ $ x_1 < x_2 <...< x_{2s+1}$
and suppose to the contrary that X is homogeneous for $C$. By
observation \ref{observ:svalues} we get that 
$\left\lfloor{\log}_s(x_1)\right\rfloor < \left\lfloor {\log}_s(x_{s+1})
\right\rfloor  < \left\lfloor {\log}_s(x_{2s+1}) \right\rfloor$ and therefore
$C(x_1,x_{s+1}) < C(x_{s+1},x_{2s+1})$ contrary to homogeneity.
\end{proof}

\bibliographystyle{plain}

\end{document}